\documentclass[twoside,12pt]{amsart}
\usepackage{amsthm,amsmath,amssymb,amscd,enumerate, ifthen, amsfonts}

\newtheorem{theorem}[subsubsection]{Theorem}
\newtheorem{corollary}[subsubsection]{Corollary}
\newtheorem{lemma}[subsubsection]{Lemma}

\newtheorem{assump}[subsubsection]{Assumption}

\newtheorem{proposition}[subsubsection]{Proposition}
\theoremstyle{definition}

\newtheorem{definition}[subsubsection]{Definition}

\theoremstyle{remark}

\theoremstyle{remark}

\newtheorem{remark}[subsubsection]{Remark}

\numberwithin{equation}{subsubsection}

\newtheorem{condition}[subsubsection]{Condition}

\newcommand{\Mbar}{\overline{\M}}
\newcommand{\proj}{\mathbb{P}}

\newcommand{\fbar}{\bar{f}}
\newcommand{\gbar}{\bar{g}}

\newcommand{\X}{\mathcal{X}}
\newcommand{\Y}{\mathcal{Y}}
\newcommand{\Z}{\mathcal{Z}}
\newcommand{\W}{\mathcal{W}}
\newcommand{\K}{\mathcal{K}}
\newcommand{\M}{\mathcal{M}}
\newcommand{\C}{\mathcal{C}}
\newcommand{\A}{\mathbb{A}}
\newcommand{\B}{\mathcal{B}}
\newcommand{\E}{\mathcal{E}}

\newcommand{\F}{\mathcal{F}}
\newcommand{\V}{\mathcal{V}}
\newcommand{\U}{\mathcal{U}}
\newcommand{\G}{\mathbb{G}}
\newcommand{\cG}{\mathcal{G}}
\newcommand{\sL}{\mathcal{L}}

\newcommand{\sH}{\mathcal{H}}
\newcommand{\sO}{\mathcal{O}}
\newcommand{\T}{\mathcal{T}}
\newcommand{\hilbf}{\underline{Hilb}}
\newcommand{\quot}{\underline{Quot}}

\def\<{\left\langle}
\def\>{\right\rangle}

\title[Space of Twisted Stable Maps to a Quotient Stack]{On The Global Quotient Structure of The Space of Twisted Stable Maps to a Quotient Stack}

\author{Dan Abramovich}
\address{Department of Mathematics\\ Brown University\\ 151 Thayer Street\\ Providence\\ RI 02912}
\email{abrmovic@math.brown.edu}
\thanks{Research of D.A. partially supported by NSF grant DMS-0335501.}

\author{Tom Graber}
\address{Department of Mathematics\\ University of California\\ Berkeley\\ CA 94720}
\curraddr{Department of Mathematics\\ California Institute of Technology\\ Mathematics 253-37\\ Caltech\\ Pasadena\\ CA 91125}
\email{graber@caltech.edu}
\thanks{Research of T.G. partially supported by NSF grant DMS-0301179 and 
an Alfred P. Sloan Research Fellowship.}

\author{Martin Olsson} 
\address{School of Mathematics\\ Institute for Advanced Study\\ 1 Einstein Drive\\ Princeton\\ NJ 08540}
\curraddr{Department of Mathematics\\ University of Texas at Austin\\ 1 University Station C1200\\ Austin\\ TX 78712-0257}
\email{molsson@math.utexas.edu}
\thanks{Research of M.O. partially supported by an NSF post--doctoral research fellowship.}

\author{Hsian-Hua Tseng}
\address{Department of Mathematics\\ University of California\\ Berkeley\\ CA 94720}
\curraddr{Department of Mathematics\\ University of British Columbia\\ 1984 Mathematics Road\\ Vancouver\\ B.C. V6T 1Z2\\
Canada}
\email{hhtseng@math.ubc.ca}

\date{\today}

\begin{document}
\begin{abstract}
Let $\X$ be a tame proper Deligne-Mumford stack of the form $[M/G]$ where $M$ is a scheme and $G$ is an algebraic group. We prove that the stack $\K_{g,n}(\X,d)$ of twisted stable maps is a quotient stack and can be embedded into a smooth Deligne-Mumford stack. When $G$ is finite, we give a more precise construction of $\K_{g,n}(\X,d)$ using Hilbert schemes and admissible $G$-covers.
\end{abstract}
\maketitle

\section{Introduction}
We fix a base scheme $B$.

Let $\X$ be a tame proper Deligne-Mumford stack of finite presentation over $B$ with projective coarse moduli space $X$. In this paper, we consider the moduli stack $\K_{g,n}(\X,d)$ of $n$-pointed genus $g$ twisted stable maps to $\X$ of degree $d$. This stack is studied in (\cite{AV}), in which it is shown to exist by verifying Artin's axioms. By (\cite{AV}, Theorem 1.4.1), $\K_{g,n}(\X,d)$ is a proper Deligne-Mumford stack with projective coarse moduli space. 

Orbifold Gromov-Witten theory (see for instance \cite{AGV}) may be interpreted as a ``virtual" intersection theory on the stack $\K_{g,n}(\X,d)$. Certain aspects of orbifold Gromov-Witten theory require a more detailed understanding of the stack $\K_{g,n}(\X,d)$. This is of interest in particular when $\X$ is a quotient stack. More precisely, we consider a stack $\X$ which satisfies the following assumption.
\begin{assump}\label{glbquot}
The stack $\X$ over $B$ is of finite presentation, tame, proper, Deligne-Mumford, with quasi-projective coarse moduli space, and is of the form $[M/G]$ where $M$ is a quasi-projective scheme and $G$ is a linear algebraic group.
\end{assump}
The following is the main result of the paper.
\begin{theorem}\label{main}
Let $\X$ be as in Assumption \ref{glbquot}. Then every connected component of $\K_{g,n}(\X,d)$ admits a locally closed immersion into a Deligne-Mumford stack of the form $[\mathcal{P}/\G]$, where $\mathcal{P}$ is a smooth quasi-projective scheme and $\G$ is an algebraic group.
\end{theorem}
Recall that a stack has the resolution property if every coherent sheaf is a quotient of a vector bundle (see for instance \cite{To}).
\begin{corollary}\label{prop}
Let $\X$ be as in Assumption \ref{glbquot}.
\begin{enumerate}

\item [{\rm (i)}] The stack $\K_{g,n}(\X,d)$ is a quotient of a scheme by an algebraic group.
\item [{\rm (ii)}] The resolution property holds for $\K_{g,n}(\X,d)$.
\end{enumerate}
\end{corollary}

We also consider similar questions in the case when $\X$ admits an action by an
algebraic group. 
\begin{assump}\label{gact}
Let $\X=[M/G]$ be as in Assumption \ref{glbquot}, and $H$ a linear algebraic group such that \begin{enumerate}
\item
there is a $G\times H$-action on $M$, which induces $H$-actions on the stack $\X$ (see for instance \cite{Ro} for a precise definition) and its coarse moduli space $X$, and 
\item
there is an $H$-linearized ample invertible sheaf on $X$.
\end{enumerate}
\end{assump}

\begin{corollary}\label{equi}
Let $\X$ be as in Assumption \ref{gact}. Then there is an $H$-equivariant locally closed immersion of $\K_{g,n}(\X,d)$ into a Deligne-Mumford stack smooth over $B$.
\end{corollary}
Corollary \ref{equi} shows that the technical assumption used in the proof of virtual localization formula (\cite{GP}) holds for $\K_{g,n}(\X,d)$. This allows one to do virtual localization calculations in orbifold Gromov-Witten theory.

When $G$ is finite, we give a more precise construction of $\K_{g,n}(\X,d)$ using Hilbert schemes and stacks of admissible $G$-covers. In this special case, Theorem \ref{main} and Corollaries \ref{prop}, \ref{equi} may be verified directly using this construction.

The paper is organized as follows. In Section \ref{pf}, we prove Theorem \ref{main} and Corollaries \ref{prop}, \ref{equi}. In Section \ref{op} we discuss some properties of the universal family of twisted stable maps that follow from our proof of Theorem \ref{main}. In Section \ref{finiquot} we give a construction of $\K_{g,n}(\X,d)$ in the case when $\X$ is a quotient by a finite group.

\section*{Acknowledgments} 
We thank R. Hartshorne, Y. Hu, A. Kresch, J. Starr and R. Vakil for helpful
conversations, and the referee for suggestions.

\section{Proof of \ref{main}, \ref{prop} and \ref{equi}}\label{pf}
\subsection{}
\begin{definition}(\cite{OS}, Section 5)
Let $\X$ be a tame Deligne-Mumford stack with coarse moduli space $\pi:\X\to X$. A locally free sheaf $\V$ on $\X$ is a {\em generator} for a quasi-coherent sheaf $\F$ if the map
$$\pi^*(\pi_*\sH om_{\sO_\X}(\V,\F))\otimes_{\sO_\X}\V\to \F$$
is surjective. A locally free sheaf $\V$ on $\X$ is a {\em generating sheaf} if it is a generator for every quasi-coherent sheaf on $\X$.
\end{definition}
\begin{remark}\label{gen}
\hfill
\begin{enumerate}
\renewcommand{\labelenumi}{(\roman{enumi})}
\item
By (\cite{OS}, Theorem 5.5), if $\X$ is a quotient of a scheme quasi-compact over $B$ by an algebraic group, then $\X$ has a generating sheaf.
\item
By (\cite{OS}, Theorem 5.2) and the Remark following it, a sheaf $\V$ is a generating sheaf of $\X$ if and only if the following holds: for each algebraically closed field $F$ and each 1-morphism $\zeta: \text{Spec}(F)\to\X$ with stabilizer group $H_\zeta$, every irreducible representation of $H_\zeta$ occurs in $\zeta^*\V$.
\end{enumerate}
\end{remark}
\begin{lemma}\label{gline} 
In the situation on Assumption \ref{gact}, there is an $H$-linearized generating sheaf of $\X$.
\end{lemma}
\begin{proof}
Choose an embedding of $G\times H$ into $GL_N$ and consider the natural map $c:\X\to BGL_N$. By (\cite{OS}, Theorem 5.5) we can find a representation $W$ of $GL_N$ such that the
pullback $c^*\W$ of the corresponding locally free sheaf $\W$ on $BGL_N$ is a generator for every quasi-coherent sheaf on $\X$. Hence $c^*\W$ is a generating sheaf of $\X$. Moreover, $H$ acts on $c^*\W$ via the induced action of $H\subset GL_N$ on $W$. Hence $c^*\W$ is a $H$-linearized generating sheaf of $\X$.
\end{proof}
\begin{lemma}\label{rep}
If $\alpha:\Y\to\X$ is a representable morphism and $\V$ is a generating sheaf of $\X$, then $\alpha^*\V$ is a generating sheaf of $\Y$.
\end{lemma}
\begin{proof}
Let $F$ be an algebraically closed field and $\zeta: \text{Spec}(F)\to\Y$ a 1-morphism with stabilizer group $G_\zeta$. By Remark \ref{gen} ({\rm ii}), it suffices to show that every irreducible representation of $G_\zeta$ occurs in $\zeta^* \alpha^*\V$. Let $G_{\alpha\circ\zeta}$ be the stabilizer group of the composite $\alpha\circ \zeta: \text{Spec}(F)\to \Y \to \X$. Since $\alpha$ is representable, the natural map $G_\zeta\to G_{\alpha\circ\zeta}$ is an inclusion of finite groups. The orders of $G_\zeta$ and $G_{\alpha\circ\zeta}$ are invertible in $F$ since $\X$ is tame. Since $\V$ is a generating sheaf of $\X$, every irreducible representation of $G_{\alpha\circ\zeta}$ occurs in $\zeta^* \alpha^*\V$. If $W$ is an irreducible representation of $G_\zeta$, consider the induced representation $Ind_{G_\zeta}^{G_{\alpha\circ\zeta}}W$ of $G_{\alpha\circ\zeta}$. Let $V$ be an irreducible representation of $G_{\alpha\circ\zeta}$ that is contained in $Ind_{G_\zeta}^{G_{\alpha\circ\zeta}}W$. By Frobenius reciprocity (see for instance \cite{FH}, Corollary 3.20), the representation $Res V$ of $G_\zeta$ obtained by restriction contains $W$. Hence every irreducible representation of $G_\zeta$ occurs in $\zeta^* \alpha^*\V$ as well. 
\end{proof}

\subsection{}\label{not}
Fix a generating sheaf $\V$ on $\X$ and an ample invertible sheaf $\sO_X(1)$ on $X$.

Let $S$ be a scheme and write $\X_S:=\X\times S$. An object $f:(\C, \{\Sigma_i\})\to \X_S$ of $\K_{g,n}(\X,d)(S)$ consists of the following data (see \cite{AV}, Definition 4.3.1):
$$ \begin{array}{ccc}
 \C             & \overset{f}{\to}&  \X_S \\
\pi\downarrow     &    & \downarrow \\

 C             & \overset{\fbar}{\to}& X_S \\
\gbar\downarrow     &    &\\
S \end{array}$$
along with $n$ closed substacks $\Sigma_i\subset \C$ such that
\begin{enumerate}
\item $\C$ is a twisted nodal $n$-pointed curve over $S$ (see \cite{AV}, Definition 4.1.2),
\item $f:\C\to \X_S$ is representable,
\item $\Sigma_i$ is an \'etale gerbe over $S$, for $i=1,...,n$, and
\item the map $\fbar: (C,\{p_i\})\to X$ between coarse moduli spaces induced from $f$ is a stable $n$-pointed map of degree $d$.
\end{enumerate}
Here $\pi: \C\to C$ is the map to the coarse moduli space, and $\gbar:C\to S$ is the structure morphism.
Let $g=\gbar\circ\pi:\C\to S$ and denote by $\omega_{C/S}$ the dualizing sheaf of $C/S$. 

Consider the invertible sheaf $\sL:=\omega_{C/S}(p_1+...+p_n)\otimes\fbar^*\sO_X(3)$. Stability of $\fbar$ implies that $\sL$ is ample relative to $S$ on $C$, see for instance (\cite{FuP}, Section 1.1). For a sheaf $\F$ on $\C$ we write $\F(s)$ for $\F\otimes \pi^*\sL^{\otimes s}$. 

\begin{lemma}\label{bdd}
There exists an integer $N$ such that for any $S$ and $f$, the following hold:
\begin{enumerate}

\item [{\rm (i)}] $R^ig_*(f^*\V^\vee(N))=0$ for all $i>0$ .

\item [{\rm (ii)}] The sheaf $g_*(f^*\V^\vee(N))$ is locally free and its formation commutes with arbitrary base change $S'\to S$.

\item [{\rm (iii)}] The composition map 
\begin{equation}\label{sur}
g^*g_*(f^*\V^\vee(N))\otimes f^*\V(-N)\to f^*(\V^\vee\otimes \V)\to \sO_\C
\end{equation}
is surjective.

\item [{\rm (iv)}] $R^i\gbar_*\sL^N=0$ for all $i>0$. The sheaf $\gbar_*\sL^N$ is locally free, its formation commutes with arbitrary base change, and the adjunction map $$\gbar^*\gbar_*\sL^N\to\sL^N$$ is surjective.

\end{enumerate}
\end{lemma}
\begin{proof}
Since $\X$ is of finite presentation over $B$, there exists a quasi-compact Noetherian scheme $B_0$ and a model $\X_0$ for $\X$ over $B_0$ such that $\X$ is obtained by a base change $B\to B_0$ (that is, $\X=\X_0\times_{B_0}B$). Therefore the stack $\K_{g,n}(\X,d)=\K_{g,n}(\X_0,d)\times_{B_0}B$ (see \cite{AV}, Proposition 5.2.1) is of finite type over $B_0$. So we may assume that the base $B$ is quasi-compact.
 
Since $\K_{g,n}(\X,d)$ is quasi-compact, there is a quasi-compact scheme $W$ and a stable map $(\C_W, \{\Sigma_i\})\to \X_W$ over $W$ such that, for any scheme $S$, the objects of $\K_{g,n}(\X,d)(S)$ are all \'etale locally on $S$ obtained by pullback via a morphism $S\to W$. In other words, there is a smooth cover $W\to \K_{g,n}(\X,d)$ with $W$ a quasi-compact scheme so that every morphism $S\to \K_{g,n}(\X,d)$ \'etale locally factors through $W\to \K_{g,n}(\X,d)$. Therefore it suffices to consider $S=W$ and a versal family $(\C, \{\Sigma_i\})\to \X_S$. In this case, ({\rm i}) follows from (\cite{H}, Theorem III.8.8), ({\rm ii}) follows from ({\rm i}) and Theorem on Cohomology and Base Change (see \cite{H}, Theorem III.12.11). ({\rm iv}) follows from (\cite{H}, Theorem III.8.8 and Theorem III.12.11). We choose $N$ large enough so that ({\rm i}), ({\rm ii}) and ({\rm iv}) all hold. 

We now prove the surjectivity of (\ref{sur}). By Lemma \ref{rep}, $f^*\V$ is a generating sheaf on $\C$. The twist $f^*\V(-N)$ is also a generating sheaf on $\C$. Hence the map
\begin{equation}\label{sur2}
\pi^*\pi_*(f^*\V^\vee(N))\otimes f^*\V(-N)\to \sO_\C
\end{equation}
is surjective. Here we use $\sH om(f^*\V(-N),\sO_\C)=f^*\V^\vee(N)$.

Since $C/S$ is projective with an ample sheaf $\sL$, we can choose a sufficiently large $N$ so that in addition to ({\rm i}), ({\rm ii}), ({\rm iv}), the map
$$\gbar^*\gbar_*\pi_*f^*\V^\vee(N)\to \pi_*f^*\V^\vee(N)$$ is surjective. Pulling back to $\C$ and tensoring with $f^*\V(-N)$, we obtain a surjective map
$$g^*g_*f^*\V^\vee(N)\otimes f^*\V(-N)\to \pi^*\pi_*f^*\V^\vee(N)\otimes f^*\V(-N).$$
Combining this with the surjection (\ref{sur2}), we obtain the surjectivity of (\ref{sur}).
\end{proof}

We assume the notation used in Section \ref{not}. The proof of Theorem \ref{main} is divided into several steps. We first rigidify the moduli problem using frames of $g_*f^*\V^\vee(N)$ and $\gbar_*\sL^{\otimes N}$. Next we prove that the original moduli stack is the stack quotient of the rigidified moduli space by an algebraic group, and the rigidified moduli space is an algebraic space. The last step is to immerse the rigidified moduli space into a scheme by a Quot functor construction.

\subsection{Rigidification}
Note that the ranks of $g_*f^*\V^\vee(N)$ and $\gbar_*\sL^{\otimes N}$ are constants on connected components of $\K_{g,n}(\X,d)$. For positive integers $l_1$ and $l_2$, let $$\K_{l_1,l_2}\subset \K_{g,n}(\X,d)$$ be the open-and-closed substack of objects for which the ranks of $g_*f^*\V^\vee(N)$ and $\gbar_*\sL^{\otimes N}$ are equal to $l_1$ and $l_2$ respectively. Define a stack $\Gamma_{l_1,l_2}$ as follows:
to any scheme $S$ we associate the groupoid of triples $(((\C,\{\Sigma_i\})\to \X_S),\iota_1,\iota_2)$ where $((\C,\{\Sigma_i\})\to \X_S)$ is an object of $\K_{g,n}(\X,d)(S)$, and $$\iota_1: \sO_S^{\oplus l_1}\to g_*f^*\V^\vee(N),\quad \iota_2: \sO_S^{\oplus l_2}\to \gbar_*\sL^{\otimes N}$$ are isomorphisms (the {\em frames} of $g_*f^*\V^\vee(N)$ and $\gbar_*\sL^{\otimes N}$). A priori $\Gamma_{l_1,l_2}$ is a 2-category, but by (\cite{AV}, Lemma 4.2.3) it is in fact equivalent to a category.
The stack $\Gamma_{l_1,l_2}$ relatively represents an Isom functor over the stack $\K_{g,n}(\X,d)$, hence it is algebraic. 
The group $GL_{l_1}\times GL_{l_2}$ acts on $\Gamma_{l_1,l_2}$ by the natural actions of $GL_{l_1}$ on $\sO_S^{\oplus l_1}$ and $GL_{l_2}$ on $\sO_S^{\oplus l_2}$. This induces an action of the diagonally embedded $\G_m\hookrightarrow GL_{l_1}\times GL_{l_2}$: an element $u\in\G_m(S)$ sends $((\C,\{\Sigma_i\})\to \X_S),\iota_1,\iota_2)$ to $((\C,\{\Sigma_i\})\to \X_S),u\cdot\iota_1,u\cdot\iota_2)$.
Let $\tilde{\K}_{l_1,l_2}:=[\Gamma_{l_1,l_2}/\G_m]$ be the stack quotient, which is algebraic by (\cite{Ro}, Theorem 4.1).
We call the stack $\tilde{\K}_{l_1,l_2}$ the rigidified moduli space.

Let $\G$ be the quotient of $GL_{l_1}\times GL_{l_2}$ by the diagonally embedded $\G_m$. The $GL_{l_1}\times GL_{l_2}$-action on $\Gamma_{l_1,l_2}$ induces an action of $\G$ on $\tilde{\K}_{l_1,l_2}$.
\begin{lemma}
$[\tilde{\K}_{l_1,l_2}/\G]\simeq \K_{l_1,l_2}$.
\end{lemma}
\begin{proof}
It is easy to see that $[\tilde{\K}_{l_1,l_2}/\G]\simeq [\Gamma_{l_1,l_2}/GL_{l_1}\times GL_{l_2}]$. There is a natural functor $$\Gamma_{l_1,l_2}\to \K_{l_1,l_2}$$ which sends a triple $(((\C,\{\Sigma_i\})\to \X_S),\iota_1,\iota_2)$ to $((\C,\{\Sigma_i\})\to \X_S)$. Two triples $(((\C,\{\Sigma_i\})\to \X_S),\iota_1,\iota_2)$ and $(((\C',\{\Sigma'_i\})\to \X_{S'}),\iota'_1,\iota'_2)$ have the same image under this functor if and only if $$((\C,\{\Sigma_i\})\to\X_S)=((\C',\{\Sigma'_i\})\to \X_{S'}).$$ In this case we have ${\iota'}^{-1}_1\circ\iota_1\in GL_{l_1}$, ${\iota'}^{-1}_2\circ\iota_2\in GL_{l_2}$. Hence $$[\Gamma_{l_1,l_2}/GL_{l_1}\times GL_{l_2}]\simeq \K_{l_1,l_2}.$$

\end{proof}

\subsection{Representability by algebraic spaces}
We show that $\tilde{\K}_{l_1,l_2}$ is represented by an algebraic space.
Consider the projective space $t:\proj_\X^{l_2-1}\to \X$ over $\X$. Clearly $\proj_\X^{l_2-1}$ is a quotient of a scheme by an algebraic group. Let $r:\A\to\X$ be the stack obtained by taking the relative spectrum of $Sym^\bullet(((t^*\V)\otimes\sO_{\proj_\X^{l_2-1}}(-1))^{\oplus l_1})$ over $\proj_\X^{l_2-1}$. The stack $\A$ is affine over $\proj_\X^{l_2-1}$, hence it is also a quotient of a scheme by an algebraic group.

\begin{lemma}\label{spaceA}
\hfill
\begin{enumerate}

\item [{\rm (i)}] If $\V$ is a generating sheaf on $\X$, then $r^*\V$ is a generating sheaf on $\A$.
\item [{\rm (ii)}] The coarse moduli space $A$ of $\A$ is affine over $\proj_X^{l_2-1}$.
\end{enumerate}
\end{lemma}
\begin{proof}
Since $r$ is representable, ({\rm i)} follows from Lemma \ref{rep}. For ({\rm ii}), note first that the coarse moduli space of $\proj_\X^{l_2-1}$ is $\proj_X^{l_2-1}$. This follows from the fact that formation of coarse moduli space commutes with flat base change on the coarse moduli space (see \cite{AV}, Lemma 2.2.2). We can work \'etale locally on $\proj_X^{l_2-1}$. Then ({\rm ii}) follows from the following fact: if $Y=\text{Spec}(R)$ is an affine $k$-scheme with the action of a finite group $G$ of order invertible in $k$, then the coarse moduli space of $[Y/G]$ is isomorphic to $\text{Spec}(R^G)$.
\end{proof}

Let $\hilbf(\A)$ be the functor over $B$ which associates to any $B$-scheme $S$ the groupoid of closed substacks $\Z\subset \A$ flat over $S$ and with proper support. Since $\A$ is a quotient of a quasi-projective scheme, by (\cite{OS}, Theorem 1.5), $\hilbf(\A)$ is representable by a disjoint union of quasi-projective schemes over $B$. 
\begin{definition}
There is a natural functor 
$$\bar{\rho}: \Gamma_{l_1,l_2}\to\hilbf(\A)$$ defined as follows: given an object $((f:(\C,\{\Sigma_i\})\to \X_S),\iota_1,\iota_2)$ of $\Gamma_{l_1,l_2}(S)$ for some scheme $S$, there is an embedding $j_C:C\to \proj_S^{l_2-1}$ given by the surjection $$\sO_C^{\oplus l_2}\overset{\iota_2}{\longrightarrow}\gbar^*\gbar_*\sL^{\otimes N}\to \sL^{\otimes N}.$$ Let $h:\C\to \proj_{\X_S}^{l_2-1}\simeq \X_S\times_S\proj_S^{l_2-1}$ be the map induced by $f$ and $j_C$. There are isomorphisms
$$(h^*(t^*\V)\otimes h^*\sO_{\proj_{\X_S}^{l_2-1}}(-1))^{\oplus l_1} \simeq \sO_\C^{\oplus l_1}\otimes f^*\V(-N)$$
$$\overset{\iota_1\otimes id}{\longrightarrow}g^*g_*(f^*\V^\vee(N))\otimes f^*\V(-N).$$
Combining this with the surjection (\ref{sur}), we obtain a closed immersion $j_\C:\C\to \A_S$ which lies over $j_C$. The functor $\bar{\rho}$ is defined by sending $((f:(\C,\{\Sigma_i\})\to \X_S),\iota_1,\iota_2)$ to $j_\C$.
\end{definition}

We analyze the fiber of this functor $\bar{\rho}$. Given $j_\C$ over a $B$-scheme $S$, note that $j_\C$ induces a morphism $C\to A$ whose composition with the projection to $\proj_S^{l_2-1}$ gives $j_C$. Pulling back the universal quotient $$\sO_{\proj_S^{l_2-1}}^{\oplus l_2}\to \sO_{\proj_{S}^{l_2-1}}(1)$$ to $C$ yields a surjection $\sO_C^{\oplus l_2}\to j_C^*\sO_{\proj_S^{l_2-1}}(1)$, where $j_C^*\sO_{\proj_S^{l_2-1}}(1)$ is isomorphic to $\sL^{\otimes N}$. Choose an isomorphism $\lambda: j_C^*\sO_{\proj_S^{l_2-1}}(1)\to \sL^{\otimes N}$. Since $C$ is connected, any two such isomorphisms differ by the multiplication by an element of the group scheme $\G_m(S)$. The composition $\sO_C^{l_2}\to j_C^*\sO_{\proj_S^{l_2-1}}(1) \overset{\lambda}{\to}\sL^{\otimes N}$ is a surjection. Applying $g_*$ yields a map $\sO_S^{\oplus l_2}\to g_*\sL^{\otimes N}$, which is the isomorphism $\iota_2$. 

The embedding $j_\C$ gives a surjection $$(h^*(t^*\V)\otimes h^*\sO_{\proj_\X^{l_2-1}}(-1))^{\oplus l_1}\to \sO_\C.$$
The isomorphisms $$(h^*(t^*\V)\otimes h^*\sO_{\proj_\X^{l_2-1}}(-1))^{\oplus l_1}\simeq \sO_\C^{\oplus l_1}\otimes f^*\V\otimes j_C^*\sO_{\proj_S^{l_2-1}}(-1)$$ and $\lambda: j_C^*\sO_{\proj_S^{l_2-1}}(1)\to \sL^{\otimes N}$ give an isomorphism $$(h^*(t^*\V)\otimes h^*\sO_{\proj_\X^{l_2-1}}(-1))^{\oplus l_1}\simeq \sO_\C^{\oplus l_1}\otimes f^*\V(-N).$$ From this we obtain a surjection $\sO_\C^{\oplus l_1}\otimes f^*\V(-N)\to \sO_\C$. This map induces a map $\sO_\C^{\oplus l_1}\to f^*\V^\vee(N)$. Applying $g_*$ yields the isomorphism $\iota_1:\sO_S^{\oplus l_1}\to g_*f^*\V^\vee(N)$. 

For $u\in\G_m(S)$, let $u\cdot$ denote the linear map given by the multiplication by $u$. If we replace $\lambda$ by $u\cdot\lambda$, then the above construction gives the pair $(u\cdot\iota_1,u\cdot\iota_2)$. 

The argument above shows that $\bar{\rho}$ factors through $\tilde{\K}_{l_1,l_2}$. Let $$\rho :\tilde{\K}_{l_1,l_2}\to \hilbf(\A)$$ be the induced functor.
For each $i$, the gerbe $\Sigma_i\subset\C$ induces a closed substack of $\A_S$ via the embedding $j_\C$. This defines a map 
$$\tau_i:\tilde{\K}_{l_1,l_2}\to\hilbf(\A).$$
Consider the map
$$\rho\times\tau_1\times...\times\tau_n: \tilde{\K}_{l_1,l_2}\to \hilbf(\A)^{n+1}.$$
This is
a morphism of algebraic stacks.
The following Lemma implies in particular that $\tilde{\K}_{l_1,l_2}$ is represented by an algebraic space (which implies that $\Gamma _{l_1, l_2}$ is also represented by an algebraic space).
\begin{lemma}\label{loc}
The map $\rho\times\tau_1\times...\times\tau_n$ is represented by locally closed immersions.
\end{lemma}

\begin{proof}
It follows from the construction that $\rho\times\tau_1\times...\times\tau_n$ exhibits $\tilde{\K}_{l_1,l_2}$ as a subfunctor of $\hilbf(\A)^{n+1}$: 
Consider an object of $\tilde{\K}_{l_1,l_2}$ represented by an object $(((\C,\{\Sigma_i\})\to \X_S),\iota_1,\iota_2)$ of $\Gamma_{l_1,l_2}$. 
An automorphism of this object corresponds to a triple $(\gamma, \delta_1, \delta_2)$ where $\gamma$ is an automorphism of the twisted stable map $((\C,\{\Sigma_i\})\to \X_S)$ and $\delta_1\in GL_{l_1}, \delta_2\in GL_{l_2}$. 
By construction, $\gamma$ yields an automorphism of $j_\C$, and $\delta_1, \delta_2$ yield an automorphism of $\A$. Furthermore, if the automorphism given by the triple $(\gamma,\delta_1,\delta_2)$ is nontrivial, so is its image under $\rho\times\tau_1\times...\times\tau_n$.

It remains to show that it is locally closed. Let $pr: \hilbf(\A)^{n+1}\to \hilbf(\A)$ be the projection to the first factor and $pr_i: \hilbf(\A)^{n+1}\to\hilbf(\A)$ the projection to the $(i+1)$-st factor, for $i=1,...,n$. Let $\Z_S/S$ be an object of $\hilbf(\A)^{n+1}$. The projection $\A\to\proj_{\X_S}^{l_2-1}\simeq \X_S\times_S\proj_S^{l_2-1}$ induces two maps $f_S:\Z_S\to \X_S$ and $j_S:\Z_S\to \proj_{S}^{l_2-1}$. Let $\fbar_S: Z_S\to X_S$ and $\bar{j}_S: Z_S\to \proj_{S}^{l_2-1}$ be the induced maps between coarse moduli spaces. The object $\Z_S/S$ is in the image of $\rho\times\tau_1\times...\times\tau_n$ if and only if
\begin{enumerate}
\item $pr(\Z_S)$ is a prestable twisted curve of genus $g$,
\item $pr_i(\Z_S)$ lie on the nonsingular locus of $pr(\Z_S)$,
\item $\fbar_S$ is a stable map of degree $d$, and
\item there is an isomorphism $$\bar{j}_S^*\sO_{\proj_{S}^{l_2-1}}(1) \longrightarrow (\omega_{Z_S/S}(\overline{pr_1(\Z_S)}+...+\overline{pr_n(\Z_S)})\otimes \fbar_S^*\sO_X(1)^{\otimes 3})^{\otimes N}$$ of line bundles, where $\overline{pr_i(\Z_S)}$ denotes the coarse moduli space of $pr_i(\Z_S)$.
\end{enumerate}
Since these are locally closed conditions, the Lemma follows.
\end{proof}

There is a natural $\G$-action on $\A$ given by the natural action of $GL_{l_2}$ on $\proj_\X^{l_2-1}$ and the fiberwise action of $GL_{l_1}$ on the bundle $\sO_{\proj_\X^{l_2-1}}(-1)^{\oplus l_1}$. The $\G$-action on $\A$ induces a $\G$-action on $\hilbf(\A)$, and hence a diagonal action on $\hilbf(\A)^{n+1}$. The map $\rho\times\tau_1\times...\times\tau_n$ is compatible with the actions of $\G$. Moreover, since $\G$ acts on $\tilde{\K}_{l_1,l_2}$ with finite stabilizers, the image of $\rho\times\tau_1\times...\times\tau_n$ is contained in the open subset of $\hilbf(\A)^{n+1}$ on which $\G$ acts with finite stabilizers.

\subsection{Representability of Quot functors}\label{emb}
We want to construct a $\G$-equivariant immersion of $\hilbf(\A)^{n+1}$ into a union of quasi-projective schemes. We do this for each component.
\begin{definition}
Let $P: K^0(\A)\to\mathbb{Z}$ be an additive homomorphism. Define $\hilbf^P(\A)$ to be the open-and-closed subscheme of $\hilbf(\A)$ classifying substacks $j:\Z\to\A$ such that the map $K^0(\A)\to \mathbb{Z}$ defined by $[E]\mapsto \chi (j^*E)$ equals $P$.
\end{definition}

\begin{lemma}
There is a $\G$-action on $\hilbf^P(\A)$ induced from that on $\hilbf(\A)$.
\end{lemma}
\begin{proof}
It suffices to show that if $\bar{k}$ is an algebraic closure of $k$, then for any element $u\in\G(\bar{k})$, $E\in K^0(\A)$, and $j: \Z\to \A_{\bar{k}}$ proper over $\bar{k}$, the integer $\chi(j^*u^*E)$ is equal to $\chi (j^*E)$.

We may assume that $k=\bar{k}$. Let $i=j\times id:\Z\times\G\to\A\times\G$. Denote by $a:\A\times \G\to\A$ the map given by the action. Put $\F:=i^*a^*E$. For any $u\in\G$, the fiber of $\F$ on $\Z$ is the sheaf $j^*u^*E$. Since Euler characteristics are constant in flat families and $\G$ is connected, the result follows.
\end{proof}

We construct for each $P$ an immersion of $\hilbf^P(\A)$ into a smooth projective scheme with compatible $\G$-action. This is done in two steps. First we immerse $\hilbf^P(\A)$ into a Quot functor over the coarse moduli space $A$, following (\cite{OS}, Section 6). Then we immerse this Quot functor into some Grassmannian, following (\cite{Gr}).

By Lemma \ref{gline}, $\proj_\X^{l_2-1}$ has a $\G$-linearized generating sheaf. Its pullback to $\A$, denoted by $\E$, is a $\G$-linearized generating sheaf of $\A$. Clearly $\proj_X^{l_2-1}$ has a $\G$-linearized ample invertible sheaf. By Lemma \ref{spaceA} (ii), $A$ is affine over $\proj_X^{l_2-1}$. Therefore a $\G$-linearized ample invertible sheaf on $\proj_X^{l_2-1}$ pulls back to a $\G$-linearized ample invertible sheaf on $A$. Choose such a sheaf $\M$ on $A$. Define a map $P_\E:K^0(\A)\to\mathbb{Z}$ by $[E]\mapsto P([\mathcal{H}om(\E,p^*E)])$, where $p:\A\to A$ is the projection. There is a natural functor
$$F_P: \hilbf^P(\A)\to \quot^{P_\E}(p_*\E^\vee/A/B)$$ defined as follows: let $\sO_\A\to J\to 0$ be an object of $$\hilbf^P(\A)=\quot^P(\sO_\A/\A/B).$$ Applying the exact functor $p_*\mathcal{H}om(\E,-)$ yields a surjection $$p_*\E^\vee\to p_*\mathcal{H}om(\E,J)$$, which is an element of $\quot^{P_\E}(p_*\E^\vee/A/B)$. The functor $F_P$ sends $\sO_\A\to J\to 0$ to $p_*\E^\vee\to p_*\mathcal{H}om(\E,J)$.

By (\cite{OS}, Proposition 6.2), this functor is an immersion. The $\G$-action on $\E$, which is compatible with that on $\A$, induces a $\G$-action on $\quot^{P_\E}(p_*\E^\vee/A/B)$. The immersion $F_P$ is $\G$-equivariant with respect to this action. 

By the argument of (\cite{Gr}, Th\'eor\`eme IV.3.2), we have the following result.
\begin{lemma}
There is an integer $M$ such that for any quotient $p_*\E^\vee\to \cG$ defining an element of $\quot^{P_\E}(p_*\E^\vee/A/B)(S)$ for some scheme $S$, the following holds (let $h_s: A_S\to S$ be the structure morphism):
\begin{enumerate}

\item [{\rm (i)}] The sheaf $h_{s*}p_*\E^\vee\otimes \M^M$ is locally free. Its formation commutes with arbitrary base change, and the map 
$$h_s^*h_{s*}p_*\E^\vee\otimes\M^{\otimes M}\to p_*\E^\vee\otimes \M^{\otimes M}$$ is surjective.
\item [{\rm (ii)}] The sheaf $h_{s*}\cG\otimes\M^{\otimes M}$ is locally free. Its formation commutes with arbitrary base change, and the map $h_s^*h_{s*}\cG\otimes\M^{\otimes M}\to \cG\otimes\M^{\otimes M}$ is surjective.
\item [{\rm (iii)}] The map $h_{s*}p_*\E^\vee\otimes\M^{\otimes M}\to h_{s*}\cG\otimes\M^{\otimes M}$ is surjective.
\item [{\rm (iv)}] The rank of $h_{s*}\cG\otimes \M^{\otimes M}$ is equal to a fixed integer $w$ independent of $\cG$.
\end{enumerate}
\end{lemma}
This Lemma implies that there is a map $$\epsilon_P:\quot^{P_\E}(p_*\E^\vee/A/B)\to Gr(p_*\E^\vee\otimes\M^{\otimes M}, w)$$ defined by $$\cG\mapsto (h_{s*}p_*\E^\vee\otimes\M^{\otimes M}\to h_{s*}\cG\otimes\M^{\otimes M}),$$where $Gr(p_*\E^\vee\otimes\M^{\otimes M}, w)$ denotes the Grassmannian of rank $w$ quotients of $p_*\E^\vee\otimes\M^{\otimes M}$.
Grothendieck's argument in (\cite{Gr}, Th\'eor\`eme IV.3.2) also shows that this map is an immersion. The $\G$-actions on $A$, $p_*\E^\vee$, and $\M$ induce a $\G$-action on $Gr(p_*\E^\vee\otimes\M^{\otimes M}, w)$, and the map $\epsilon_P$ is $\G$-equivariant. The composite map $\epsilon_P\circ F_P$ is a $\G$-equivariant immersion of $\hilbf^P(\A)$ into a Grassmannian, which is a smooth scheme as desired.

\subsection{Conclusion of the proof of \ref{main}}\label{conclusion}
By Lemma \ref{loc} there is a closed immersion $\rho\times\tau_1\times...\times\tau_n: \tilde{\K}_{l_1,l_2}\to \hilbf(\A)^{n+1}$. The image of $\rho$ is contained in a connected component of $\hilbf(\A)$ with Hilbert polynomial $P$. For each $i$ the image of $\tau_i$ is contained in a connected component of $\hilbf(\A)$ with Hilbert polynomial $P_i$. Therefore we have an immersion, 
$$\rho\times\tau_1\times...\times\tau_n: \tilde{\K}_{l_1,l_2}\to \hilbf^P(\A)\times \hilbf^{P_1}(\A)\times...\times\hilbf^{P_n}(\A)$$ which is $\G$-equivariant.
Let $\rho'=\epsilon_P\circ F_P:\hilbf^P(\A)\to Gr$ and $\tau'_i=\epsilon_{P_i}\circ F_{P_i}:\hilbf^{P_i}(\A)\to Gr_i$ be the immersions into Grassmannians constructed in Section \ref{emb}. The composition
$$(\rho'\times\tau'_1\times...\times\tau'_n)\circ(\rho\times\tau_1\times...\times\tau_n):\tilde{\K}_{l_1,l_2}\to Gr\times Gr_1\times...\times Gr_n$$ 
is a $\G$-equivariant locally closed immersion, and the image of $\tilde{\K}_{l_1,l_2}$ is contained in the locus $\mathcal{P}$ on which $\G$-acts with finite stabilizers. The scheme $\mathcal{P}$ is smooth and quasi-projective. This completes the proof of Theorem \ref{main}.

\subsection{Proof of \ref{prop}}
For the first statement of Corollary \ref{prop}, let $\K\subset \K_{g,n}(\X,d)$ be a connected component and $\K\to [\mathcal{P}/\G]$ be the immersion as in Theorem \ref{main}. The product $K:=\K\times_{[\mathcal{P}/\G]}\mathcal{P}$ is a scheme, and $\K\simeq [K/\G]$.

For the second statement, it suffices to show that any $\G$-equivariant coherent sheaf on $K$ is a quotient of a $\G$-equivariant vector bundle. Fix a $\G$-equivariant ample line bundle $L$ on $K$. Let $V$ be a $\G$-equivariant coherent sheaf of $K$. For a sufficiently large integer $k$ the sheaf $V\otimes L^{\otimes k}$ is generated by $\G$-equivariant global sections. This gives a $\G$-equivariant surjection
$$H^0(K,V\otimes L^{\otimes k})\otimes \sO_K \to V\otimes L^{\otimes k},$$
which induces a $\G$-equivariant surjection $H^0(K,V\otimes L^{\otimes k})\otimes L^{\otimes -k} \to V$. The statement follows.

\subsection{Proof of \ref{equi}}
By Lemma \ref{gline}, we can choose $\V$ and $\E$ to be $H$-equivariant. The existence of $H$-linearized ample invertible sheaf on $X$ implies that $\sL$ and $\M$ can be chosen to be $H$-equivariant. Then the immersion $\tilde{\K}_{l_1,l_2}\to Gr\times Gr_1\times...\times Gr_n$ is $H\times\G$-equivariant. This proves Corollary \ref{equi}.

\section{Other Properties}\label{op}
\subsection{}

Consider the universal family of twisted stable maps to $\X$:
\begin{equation}\label{unifamily}
\begin{CD}
  \U_{g,n}(\X,d) @>e>>\X \\
  @V \phi VV\\
 \K_{g,n}(\X,d).
\end{CD}
\end{equation}
Let $\U_{l_1,l_2}\to\K_{l_1,l_2}$ the pullback family over $\K_{l_1,l_2}$. Our construction allows us to prove another property of $\U_{l_1,l_2}\to\K_{l_1,l_2}$.
\begin{proposition}\label{embb}
Let $\X=[M/G]$ be as in Assumption \ref{glbquot}. In addition, assume $\X$ is smooth. Then there is a nonsingular Deligne-Mumford stack $\sH$ with a proper and flat family $\U_\sH\to\sH$ and an embedding $\K_{l_1,l_2}\to\sH$ such that
\begin{enumerate}

\item [{\rm (i)}] The restriction of $\U_\sH\to\sH$ to $\K_{l_1,l_2}$ is $\U_{l_1,l_2}\to\K_{l_1,l_2}$.

\item [{\rm (ii)}] The Kodaira-Spencer map (see for instance \cite{Vi}, Definition 6.9) $T_m\sH\to Ext^1(\Omega_{{\U_\sH}_m},\sO_{{\U_\sH}_m})$ is surjective for all $m\in \sH$.

\end{enumerate}
\end{proposition}
\begin{proof}
 
Our argument is parallel to that of \cite{FaP}. We begin by considering the following situation: Let $\Y$ be a smooth Deligne-Mumford stack of the form $[Y/\mathbf{G}]$ with $Y$ a quasi-projective variety and $\mathbf{G}$ a linear algebraic group. Let $$\C_S\subset \Y\times S\to S$$ be a projective flat family of twisted curves over a quasi-projective base scheme $S$. Choose a $\mathbf{G}$-linearized very ample line bundle $\sO_Y(1)$ on $Y$ and let $\sL$ denote the corresponding line bundle on $\Y$. The line bundle $\sO_Y(1)$ induces a $\mathbf{G}$-equivariant embedding $Y\hookrightarrow \proj^r$, which descends to an embedding $\Y\hookrightarrow [\proj^r/\mathbf{G}]$.

For $s\in S$ let $\C_s$ denote the fiber of $\C_S\to S$ over $s$. By (\cite{H}, Theorem III.8.8), there exists an integer $m$ such that $$H^1(\C_s,\sL^m)=0$$ for all $s\in S$. Put $N=\binom{r+m}{m}-1$. The Veronese embedding $\proj^r\to \proj^N$ given by the line bundle $\sO_{\proj^r}(m)$ is $\mathbf{G}$-equivariant, hence descends to an embedding $[\proj^r/\mathbf{G}]\hookrightarrow [\proj^N/\mathbf{G}]$. Let $W\subset \proj^N$ denote the open subscheme on which $\mathbf{G}$ acts with finite stabilizers. The stack $[W/\mathbf{G}]$ is Deligne-Mumford. The composition $$\Y\hookrightarrow [\proj^r/\mathbf{G}]\hookrightarrow [\proj^N/\mathbf{G}]$$ has image in $[W/\mathbf{G}]$. This yields a morphism $S\to \hilbf([W/\mathbf{G}])$.
\begin{lemma}
The image of $S\to \hilbf([W/\mathbf{G}])$ is contained in the open subscheme $\sH^o$ where $\hilbf([W/\mathbf{G}])$ is smooth over the base $B$.
\end{lemma}
\begin{proof}
Let $\C$ be a fiber of $\C_S\to S$. Since $\C$ is an LCI stack and $[W/\mathbf{G}]$ is smooth, the morphism $f:\C\to [W/\mathbf{G}]$ is LCI. The Lemma follows if we can show that the obstruction for $f$ vanishes. Let $N_{\C/[W/\mathbf{G}]}$ be the normal bundle of $\C\subset [W/\mathbf{G}]$ and $T_{[W/\mathbf{G}]}$ the tangent bundle of $[W/\mathbf{G}]$. It follows from (\cite{Il} Proposition 2.2.2 and Corollarie 3.2.7) that the obstruction lies in $Ext^1(N_{\C/[W/\mathbf{G}]}^\vee,\sO_\C)$. Hence we need to show that 
\begin{equation}\label{van}
Ext^1(N_{\C/[W/\mathbf{G}]}^\vee,\sO_\C)=0.
\end{equation}
Consider the exact sequence $$0\to N_{\C/[W/\mathbf{G}]}^\vee\to \Omega^1_{[W/\mathbf{G}]}|_\C\to \Omega^1_\C\to 0.$$ Applying $Hom(-,\sO_\C)$ and taking cohomology yield $$Ext^1(\Omega^1_{[W/\mathbf{G}]},\sO_\C)\to Ext^1(N_{\C/[W/\mathbf{G}]}^\vee,\sO_\C)\to Ext^2(\Omega^1_\C,\sO_\C).$$
By (\cite{ACV}, Section 3.0.3), $Ext^2(\Omega^1_\C,\sO_\C)=0$. Also, $$Ext^1(\Omega^1_{[W/\mathbf{G}]},\sO_\C)=Ext^1(\sO_\C,T_{[W/\mathbf{G}]})=H^1(\C, f^*T_{[W/\mathbf{G}]}).$$ Therefore (\ref{van}) follows from $$H^1(\C, f^*T_{[W/\mathbf{G}]})=0.$$ 
Since $\pi: W\to [W/\mathbf{G}]$ is smooth, there is a surjection $T_W\to \pi^* T_{[W/\mathbf{G}]}$. This descends to a surjection $\T\to T_{[W/\mathbf{G}]}$ on $[W/\mathbf{G}]$, yielding an exact sequence of locally free sheaves $$0\to Ker \to \T \to T_{[W/\mathbf{G}]}\to 0.$$
Pulling back to $\C$ and taking cohomology yield an exact sequence $$H^1(\C, f^*\T)\to H^1(\C,f^*T_{[W/\mathbf{G}]})\to H^2(\C, Ker).$$
On the other hand, the Euler sequence $0\to \sO_{\proj^N}\to \sO_{\proj^N}(1)^{\oplus N+1}\to T_{\proj^N}\to 0$ gives an exact sequence of locally free sheaves on $[W/\mathbf{G}]$,
$$0\to \sO_{[W/\mathbf{G}]}\to \sO_{[W/\mathbf{G}]}(1)^{\oplus N+1}\to \T\to 0.$$
Pulling back to $\C$ and taking cohomology yield an exact sequence $$H^1(\C, \sL^m)^{\oplus N+1}\to H^1(\C,f^*\T)\to H^2(\C, f^*\sO_{[W/\mathbf{G}]}).$$
Since $\C$ is tame and one-dimensional, both $H^2(\C,Ker)$ and $H^2(\C,f^*\sO_{[W/\mathbf{G}]})$ vanish. By the choice of $m$, $H^1(\C,\sL^m)=0$. It follows that $H^1(\C, f^*T_{[W/\mathbf{G}]})=0$ and the Lemma is proved.
\end{proof}
Let $\U^o\to \sH^o$ be the universal family of the Hilbert functor and $S\hookrightarrow Z$ a closed embedding in a nonsingular scheme $Z$. The vanishing (\ref{van}) implies that the family $\U^o\times Z\to \sH^o\times Z$ satisfied (ii) of this Proposition, and the restriction of $\U^o\times Z\to \sH^o\times Z$ to $S$ is $\C_S\to S$.

In Section \ref{pf}, the stack $\K_{l_1,l_2}$ is constructed as a quotient by $\G$ of a (locally closed) subscheme $\tilde{\K}_{l_1,l_2}$ of a Hilbert scheme. The universal family $\U_{l_1,l_2}\to \K_{l_1,l_2}$ is the stack quotient of the universal family $\tilde{\U}_{l_1,l_2}\subset \A\times \tilde{\K}_{l_1,l_2}\to \tilde{\K}_{l_1,l_2}$ of the Hilbert scheme. By the definition of $\A$, the $\G$-action on $\A$ is a fiberwise action with respect to the morphism $r: \A\to\X$. It follows that we may choose $\sL$ to be $\G$-linearized. The nonsingular embedding of $\tilde{\K}_{l_1,l_2}$ constructed in Section \ref{conclusion} is $\G$-equivariant. We conclude that the above construction can be applied $\G$-equivariantly to produce the required family $\U_\sH\to \sH$.
\end{proof}

\begin{remark}
By construction, the family $\U_\sH\to \sH$ can be embedded into a projective bundle.
\end{remark}

\begin{corollary}
In the situation of Proposition \ref{embb}, the image in $\sH$ of the nodal locus of $\U_\sH$ is a normal crossing divisor.
\end{corollary}
\begin{proof}
This follows from Proposition \ref{embb} (ii), deformation theory of twisted curves (see \cite{AV}, \cite{ACV}), and the arguments of (\cite{DM}, Section 1).
\end{proof}

\subsection{Another Moduli Stack}
In (\cite{AGV}), another moduli stack is considered. Let $\Mbar_{g,n}(\X,d)$ be the moduli stack of twisted stable maps to $\X$ with sections to all gerbes. In (\cite{AGV}), this moduli stack is used to give a more transparent construction of orbifold Gromov-Witten invariants. The arguments given in Section \ref{pf} and \ref{op} can be easily modified to show that our main results, Theorem \ref{main}, Corollaries \ref{prop}, \ref{equi}, and Proposition \ref{embb} hold for $\Mbar_{g,n}(\X,d)$.

\section{Global quotients by finite groups}\label{finiquot}
In this section, we consider stacks which are global quotients of schemes by finite groups. Let $\X$ be a tame proper Deligne-Mumford stack of the form $[M/G]$, where $M$ is a projective variety and $G$ is a finite group acting on $M$. The coarse moduli space of $\X$ is denoted by $M/G$ and is assumed to be projective. The goal of this section is to give a construction of the stack $\K_{g,n}(\X,\beta)$ using projective geometry in this case.

Our construction is based on the following observation. Let $\C\to [M/G]$ be a twisted stable map. Consider the Cartesian diagram,
$$ \begin{array}{ccc}
  \C\times_{[M/G]}M & \to   & M\\
\downarrow     &    &\downarrow \\
 \C             & \overset{f}{\to}& [M/G] \end{array}$$
where $M\to [M/G]$ is the atlas. Since $\C\to [M/G]$ is representable, the fiber product $\C\times_{[M/G]}M$ is a scheme. The composition of the map $\C\times_{[M/G]}M\to \C$ and the map $\C\to C$ to the coarse moduli space gives a map $\C\times_{[M/G]}M\to C$, which is an admissible $G$-cover (see \cite{ACV}). The map $\C\times_{[M/G]}M\to M$ is $G$-equivariant. Also, there is a induced (stable) map $C\to M/G$ between coarse moduli spaces. To parametrize twisted stable maps $\C\to [M/G]$, it suffices to parametrize three items:
\begin{enumerate}
\item
stable maps $C \to M/G$ to the coarse moduli space $M/G$,
\item 
admissible $G$-covers over the curves $C$,
\item morphisms from the domain of an admissible $G$-cover to $M$ which are $G$-equivariant.
\end{enumerate}
Since the curves $C$ can vary, we need to parametrize admissible $G$-covers to all possible target curves $C$.
\subsection{}
Let $\mathfrak{M}_{g,n}$ be the Artin stack of prestable $n$-pointed genus $g$ curves (see for instance \cite{Be}), and $Adm_{g,n}(G)$ the stack of admissible $G$-covers over prestable $n$-pointed genus $g$ curves. It follows from the construction of (\cite{ACV}, Appendix B) that $Adm_{g,n}(G)$ is an Artin stack. There is a morphism $$F_1:Adm_{g,n}(G)\to\mathfrak{M}_{g,n}$$ given by sending an admissible $G$-cover $E\to C$ to the curve $C$. Denote by $\K_{g,n}(M/G,\beta)$ the stack of $n$-pointed genus $g$ stable maps of degree $\beta$ to $M/G$. There is another morphism $$F_2:\K_{g,n}(M/G, \beta)\to\mathfrak{M}_{g,n}$$ given by sending a stable map to its domain curve. Let $$\C_{g,n}^G\to Adm_{g,n}(G),\quad \T_{g,n}^G\to Adm_{g,n}(G)$$ denote the universal domain and the universal target respectively. There is a diagram, 
$$ \begin{array}{ccc}
 \C_{g,n}^G             & \to&  \T_{g,n}^G \\
\downarrow     &  \swarrow  &  \\
 Adm_{g,n}(G)            & &  \end{array}$$
which is the universal admissible $G$-cover: for a morphism $S\to Adm_{g,n}(G)$, the pullback of this diagram via this morphism gives a diagram, 
$$ \begin{array}{ccc}
 E            & \to&  C \\
\downarrow     &  \swarrow  &  \\
 S            & &  \end{array}$$
which is an $S$-family of admissible $G$-covers.

Let $F_3$ denote the composition $\C_{g,n}^G\to Adm_{g,n}(G)\overset{F_1}{\to}\mathfrak{M}_{g,n}$. Put $$\B:=Adm_{g,n}(G)\times_{F_1, \mathfrak{M}_{g,n}, F_2}\K_{g,n}(M/G,\beta)$$ and let $$\E:=\C_{g,n}^G\times_{F_3,\mathfrak{M}_{g,n},F_2}\K_{g,n}(M/G,\beta)\to\B$$ be the pullback family.
\begin{lemma}
\hfill
\begin{enumerate}
\item [{\rm (i)}]
The stack $\B$ is Deligne-Mumford. 
\item [{\rm (ii)}]
The family $\E\to\B$ is projective.
\end{enumerate}
\end{lemma}
\begin{proof}
For {\rm (i)}, since $\K_{g,n}(M/G,\beta)$ is Deligne-Mumford, it suffices to show that the morphism $Adm_{g,n}(G)\to\mathfrak{M}_{g,n}$ is of Deligne-Mumford type (see \cite{LMB}, Remarque 7.3.3). Let $E\to C\to S$ be an admissible $G$-cover over an $n$-pointed genus $g$ nodal curve $C\to S$. Let $\C=[E/G]$ be the stack quotient of $E$ by $G$. Then $\C\to S$ is an $n$-pointed genus $g$ twisted curve and $E\to S$ induces a morphism $\C\to BG$. By (\cite{ACV}, Theorem 4.3.2), the sheaf of automorphisms of $E\to S$ fixing $C\to S$ is the same as the sheaf of automorphisms of the morphism $\C\to BG$ fixing $C\to S$. An infinitesimal automorphisms of $\C\to BG$ uniquely gives an infinitesimal automorphisms of $\C\to S$. It is shown in (\cite{ACV}, Section 3.0.4) that the sheaf of infinitesimal automorphisms of $\C\to S$ is the same as the sheaf of infinitesimal automorphisms of $C\to S$. It follows that there is no infinitesimal automorphisms of $\C\to BG$ fixing $C\to S$, and $Adm_{g,n}(G)\to\mathfrak{M}_{g,n}$ is of Deligne-Mumford type (\cite{LMB}, Th\'eor\`eme 8.1).

For {\rm (ii)}, we need to find a line bundle on $\E$ that is ample relative to $\B$. Let $V$ be an ample line bundle on $M/G$. For a scheme $S$, consider a stable map $f:(C, x_1,...,x_n)\to M/G$ over $S$ and an $S$-family of admissible $G$-covers $\phi :E\to C$. Stability implies that $L=\omega_{C/S}(x_1+...+x_n)\otimes f^*V^{\otimes 3}$ is relatively ample. Since $\phi$ is a finite morphism, $\phi^*L$ is relatively ample. The line bundle on $\E$ defined by $\phi^*L$ therefore is relative ample.
\end{proof}
Consider the Hom-stack $Hom^d_{\B}(\E, M\times\B)$, which is a stack over $\B$ defined as follows: for a morphism $T\to\B$ from a scheme $T$, the objects of $Hom^d_{\B}(\E, M\times\B)(T)$ are commutative diagrams 
$$ \begin{array}{ccc}
 \E\times_\B T            & \to&  M\times T \\
\downarrow     &  \swarrow  &  \\
 T            & &  \end{array}$$
such that for each $t\in T$ the morphism $\E\times_\B T\to M\times T$ restricted to the fiber over $t$ has degree $d=p^*\beta$. This Hom-stack can be constructed by descent techniques and the Hom-schemes for schemes (\cite{Gr}, IV. 4). The group $G$ acts on this Hom-stack in the following way: for a map $[f]\in Hom^d_{\B}(\E, M\times\B)$ and $\alpha\in G$, $\alpha\cdot [f]=[\alpha^{-1}\cdot f(\alpha\cdot)]$. Let $Hom^d_{\B}(\E, M\times\B)^G$ be the substack consisting of $G$-fixed points. Since $G$-equivariance is a closed condition, this is a closed substack. Objects of $Hom^d_{\B}(\E, M\times\B)^G$ are triples
$$( E\to C,\, E\to M,\, C\to M/G),$$
where
\begin{enumerate}
\item $E \to C$ is an admissible $G$-cover,
\item $E\to M$ is a $G$-equivariant morphism, and
\item $C \to M/G$ is a stable map.
\end{enumerate}
Inside $Hom^d_{\B}(\E, M\times\B)^G$ there is a substack whose objects are diagrams 
such that the map $C=E/G\to M/G$ induced from $E\to M$ coincides with the given stable map $C\to M/G$. This substack is $\K_{g,n}([M/G], \beta)$. Note that since the condition that the map $C=E/G\to M/G$ coincides with the given stable map $C\to M/G$ is a closed condition, and the stack $\K_{g,n}(M/G, \beta)$ is proper, $\K_{g,n}([M/G],\beta)$ is a closed substack of $Hom^d_{\B}(\E, M\times\B)^G$.
\subsection{}
The base $\B$ parametrizes objects
$$ \begin{array}{ccc}
 E             \\
\downarrow     &    & \\
 (C,x_1,...,x_n)             & \overset{f}{\to}& M/G \end{array}$$
 where
\begin{enumerate}
\item $E \to C$ is an admissible $G$-cover, and
\item $f: (C,x_1,...,x_n) \to M/G$ is a stable map in $\K_{g,n}(M/G,\beta)$.
\end{enumerate}

By (\cite{ACV}, Appendix B.1), the stack $Adm_{g,n}(G)$ breaks into union of open-and-closed substacks $Adm_{g,n}^{\mathbf{d}, \mathbf{e}}(G)$ which parametrizes admissible $G$-covers of type $(g,n,\mathbf{d}=(d_1,...,d_n), \mathbf{e}=(e_1,...,e_n))$. Put $\B_{\mathbf{d},\mathbf{e}}:=Adm_{g,n}^{\mathbf{d}, \mathbf{e}}(G)\times_{\mathfrak{M}_{g,n}}\K_{g,n}(M/G,\beta)\subset\B$ and denote by $\E_{\mathbf{d},\mathbf{e}}$ the pullback of $\E$ to $\B_{\mathbf{d},\mathbf{e}}$. Accordingly the stack $\K_{g,n}([M/G],\beta)$ breaks into a union of substacks $\K_{g,n}^{\mathbf{d}, \mathbf{e}}([M/G],\beta)$. We give an explicit construction of $\K_{g,n}^{\mathbf{d}, \mathbf{e}}([M/G],\beta)$. 
\subsubsection{}\label{constable}
We review the construction of ordinary stable map spaces $\K_{g,n}(Y,\beta)$ given in (\cite{FuP}). (Note that \cite{FuP} this stack is denoted by $\Mbar_{g,n}(Y,\beta)$.) Let $Y$ be a projective variety over the base $B$ and $V$ an ample line bundle on $Y$. Let $f:(C, x_1,...,x_n)\to Y$ be a $n$-pointed stable map of genus $g$ and degree $\beta$ to $Y$. Stability implies that $L=\omega_C(x_1+...+x_n)\otimes f^*V^{\otimes 3}$ is ample on $C$. An argument similar to that in the proof of (\cite{DM}, Theorem 1.2) shows that there is an integer $r>0$, which depends on $g,n, \beta$, and $Y$ but not on $C$ and $f$, such that $L^{\otimes r}$ is very ample and $h^1(C, L^{\otimes r})=0$. The line bundle $L^{\otimes r}$ gives an embedding of $C$ into a projective space $\mathbb{P}=\mathbb{P}_B^{l-1}$ where $l=h^0(C,L^{\otimes r})$. Let $Hilb(\proj\times Y)$ be the Hilbert scheme of genus $g$ curves in $\mathbb{P}\times Y$ of bidegree $(r(2g-2+n+3V\cdot\beta), \beta)$. For each $i$, let $P_i$ be the Hilbert scheme of a point in $\mathbb{P}\times Y$. Let $$I\subset H=Hilb(\proj\times Y)\times P_1\times...\times P_n$$ be the closed subscheme corresponding to the locus where the $n$ points lie on the curve. Let $J\subset I$ be the subscheme corresponding to pairs $(C,\{p_i\})$ satisfying the following conditions:
\begin{condition}\label{stablemap}
\hfill
\begin{enumerate}
\renewcommand{\labelenumi}{(\roman{enumi})}
\item The curve $C$ is prestable. 
\item The natural projection $C\to\mathbb{P}$ is a non-degenerate embedding.
\item The $n$ points $\{p_i\}$ lie on the nonsingular locus of $C$.
\item The line bundle $\mathcal{O}_{\mathbb{P}}(1)\otimes V|_C$ is isomorphic to $\omega_C(p_1+...+p_n)^{\otimes r}\otimes V^{\otimes 3r+1}|_C$.
\end{enumerate}
\end{condition}
The group $\mathbf{G}=PGL(l)$ acts on $H$ via its action on $\mathbb{P}$. The stack $\K_{g,n}(Y,\beta)$ is the quotient $[J/\mathbf{G}]$.

\subsubsection{}
Assume the notations in Section \ref{constable}. The construction of (\cite{FuP}) discussed above applied to the case $Y=M/G$ allows us to express $\K_{g,n}(M/G,\beta)$ as $[J/\mathbf{G}]$. We now construct the stack $\K_{g,n}^{\mathbf{d}, \mathbf{e}}([M/G],\beta)$. The idea is as follows: a point on $J$ corresponds to a stable map $C\to M/G$ with an embedding $C\to \proj$. Admissible $G$-covers over $C$ may be viewed as stable maps to $\mathbb{P}$. Taking quotient by $\mathbf{G}$ removes the choices of the embedding $C\to \proj$.

Let $\Mbar(\mathbb{P})$ be the stack of stable maps to $\mathbb{P}$ of the following type: 
\begin{enumerate}
\item degree $r\cdot(2g-2+n+3V\cdot\beta)$,
\item genus $h=1+|G|(g-1)+(1/2)\sum_{i=1}^n d_i(e_i-1)$,
\item the number of marked points is equal to $\tilde{n}:=\sum_{i=1}^n d_i$.
\end{enumerate}
We apply the construction in (\cite{FuP}) discussed above to $Y=\proj$. There is another projective space $\tilde{\proj}$ such that for any stable map $f: C'\to\mathbb{P}$ in
$\Mbar (\proj)$, the source curve $C'$ is embedded in $\tilde{\proj}$. Let $Hilb(\tilde{\proj}\times \proj)$ be the Hilbert scheme of genus $h$ curves in $\tilde{\proj}\times
\proj$ and, for each $i=1,...,\tilde{n}$, let $\tilde{P}_i$ be the Hilbert scheme of a point in $\tilde{\proj}\times\proj$. There is a subscheme $$J'\subset
Hilb(\tilde{\proj}\times\proj)\times \tilde{P}_1\times...\times\tilde{P}_{\tilde{n}},$$ defined by a variant of Condition \ref{stablemap}. More precisely, the subscheme $J'$
corresponds to pairs $(C', \{\tilde{p}_i\})$ satisfying the following conditions: (i) $C'$ is prestable; (ii) the map $C'\to \tilde{\proj}$ is a non-degenerate embedding; (iii)
the $\tilde{n}$ points $\{\tilde{p}_i\}$ lie on the nonsingular locus of $C'$; (iv) the line bundle $\sO_{\tilde{\proj}}(1)\otimes \sO_\proj(1)|_{C'}$ is isomorphic to
$\omega_{C'}(\tilde{p}_1+...+\tilde{p}_{\tilde{n}})^{\otimes r}\otimes \sO_\proj(1)^{\otimes 3r+1}|_{C'}$. And we have $\Mbar(\mathbb{P})=[J'/\mathbf{G}']$. Inside $\Mbar(\mathbb{P})$ there is the locus which consists of maps that are admissible $G$-covers of type $(g,n,\mathbf{d}, \mathbf{e})$. This locus can be written as $[J''/\mathbf{G}']$ for some subscheme $J''\subset J'$. A point in $J''$ gives a stable map $f: C'\to\mathbb{P}$ such that $C'\to f(C')$ is an admissible $G$-cover of type $(g,n,\mathbf{d}, \mathbf{e})$. The universal family of the Hilbert scheme $Hilb(\tilde{\proj}\times \proj)$ yields a family $U''\to J''$.

Let $Hilb(\proj)$ be the Hilbert scheme of curves of genus $g$ and degree $r\cdot(2g-2+n+3V\cdot\beta)$. There is a morphism $$J''\to Hilb(\proj)$$ given by sending $(C'\to f(C'))$ to $f(C')$. There is a morphism $$J\to Hilb(\proj)$$ defined as follows: The subscheme $J$ corresponds to pairs $(C, \{p_i\})$ satisfying Condition \ref{stablemap}. Forgetting the points $\{p_i\}$ from the pair $(C, \{p_i\})$, we obtain a curve $C$. Therefore there is a morphism $J\to Hilb(\proj\times M/G)$, which corresponds to a family $C_J\to J$ of curves in $\proj\times M/G$. For each point $s\in J$, the fiber of $C_J\to J$ over $s$ is an embedded curve in $\proj$. Hence the projection $\proj\times M/G\to \proj$ yields a $J$-family of curves in $\proj$. This defines a morphism $J\to Hilb(\proj)$.

It now follows that $\B_{\mathbf{d},\mathbf{e}}$ is the quotient $[J''\times_{Hilb(\proj)}J/\mathbf{G}'\times\mathbf{G}]$ and the universal domain $\E_{\mathbf{d},\mathbf{e}}$ is the quotient
$[U''\times_{Hilb(\proj)}J/\mathbf{G}'\times\mathbf{G}]$.

Now the stack $Hom^d_{\B_{\mathbf{d},\mathbf{e}}}(\E, M\times\B_{\mathbf{d},\mathbf{e}})$ is the quotient 
$$[Hom_{J''\times_{Hilb(\proj)}J}^d(U''\times_{Hilb(\proj)}J, M\times J''\times_{Hilb(\proj)}J)/\mathbf{G}'\times\mathbf{G}],$$
and the stack $Hom^d_{\B_{\mathbf{d},\mathbf{e}}}(\E, M\times\B_{\mathbf{d},\mathbf{e}})^G$ is the quotient 
$$[Hom_{J''\times_{Hilb(\proj)}J}^d(U''\times_{Hilb(\proj)}J, M\times J''\times_{Hilb(\proj)}J)^G/\mathbf{G}'\times\mathbf{G}].$$
The stack $\K_{g,n}^{\mathbf{d}, \mathbf{e}}([M/G],\beta)$ sits in $Hom^d_{\B_{\mathbf{d},\mathbf{e}}}(\E, M\times\B_{\mathbf{d},\mathbf{e}})^G$ in the way described above. In particular, $\K_{g,n}^{\mathbf{d}, \mathbf{e}}([M/G],\beta)$ is a quotient stack.

Theorem \ref{main}, and Corollaries \ref{prop}, \ref{equi} in the case considered in this section directly follow from this construction.

\end{document}